\documentclass[fleqn,a4paper, 12pt]{article}

\def\beq{\begin{equation}\label}
\def\eeq{\end{equation}}
\def\prend{\hfill$\blacksquare$}
\newcommand{\FB}{{\sf FB}}
\newcommand{\LAS}{{\sf LAS}}
\newcommand{\LAST}{{\sf LAST}}
\newcommand{\FIFO}{{\sf FIFO}}
\newcommand{\LIFO}{{\sf LIFO}}
\newcommand{\PS}{{\sf PS}}
\newcommand{\SRPT}{{\sf SRPT}}
\newcommand{\thez}[1]{}\newcommand{\pap}[1]{{#1}}
\newtheorem{theorem}{Theorem}
\newtheorem{lemma}[theorem]{Lemma}
\newtheorem{definition}[theorem]{Definition}
\newtheorem{proposition}[theorem]{Proposition}
\newtheorem{corollary}[theorem]{Corollary}
\newtheorem{theorem1}{Theorem}
\newtheorem{assumption}[theorem1]{Assumption}
\newcounter{figures}

\usepackage{amsfonts, amsmath, amssymb, epsfig, psfrag}

\author{M. Mandjes\thanks{CWI, Amsterdam, The Netherlands, and University of Twente, Faculty of Mathematical
 Sciences, The Netherlands\newline} \and M. Nuyens\thanks{KdV Institute for Mathematics, University of Amsterdam, The Netherlands}}

\title{Sojourn times in the M/G/1 FB queue with light-tailed service times}

\begin{document}
\maketitle
\begin{abstract}
The asymptotic decay rate of the sojourn time of a customer in the stationary  M/G/1 queue under the Foreground-Background (\FB) service discipline is studied. The \FB\ discipline gives service to those customers that have received the least service so far. We prove that for light-tailed service times the decay rate of the sojourn time is equal to the decay rate of the busy period. It is shown that \FB\ minimises the decay rate in the class of work-conserving disciplines.
\end{abstract}
\section{Introduction}
The sojourn time of a customer, i.e.~the time between his
arrival and departure, is an often used performance measure
for queues. In this \thez{chapter}\pap{note} we compute  \pap{the asymptotic decay rate of}  
the tail of the sojourn-time \index{sojourn time}distribution of the stationary M/G/1 queue with the
Foreground-Background (\FB) discipline. This decay rate is then used to compare the performance of \FB\ with other service disciplines like \PS\ and \FIFO. 

The \FB\ discipline gives service to those customers who have received the least amount of service so far. If there are $n$ such customers, each of them is served at rate $1/n$. Thus, when the {\em age} of a customer is the amount of service a customer has received, the \FB\ discipline gives priority to the {\em youngest} customers. In the literature this discipline has been called \LAS\ or  \LAST\ (least-attained service time first) as well.

Let $V$ denote the sojourn time of a customer in the stationary M/G/1 \FB\ queue.
N\'{u}\~{n}ez Queija \cite{sindothes} showed that for service-time distributions with regularly varying tails of index $\eta\in (1,2)$, the distribution of $V$ satisfies 
\beq{nuntum}P(V>x)\sim P(B>(1-\rho)x), \qquad x\to\infty,\eeq 
where $\rho$ is the load of the system, $B$ is the generic service time, and
$\sim$ means that the quotient converges to 1.  Using N\'{u}\~{n}ez Queija's  
\thez{method,  in the second part of this chapter we prove}\pap{method, Nuyens \cite{nuyenstheziz} obtained} 
(\ref{nuntum}) under weaker assumptions. In case of regularly varying service times the tail of $V$ under other disciplines, like \FIFO, \LIFO, \PS\ and \SRPT, has been found to be  heavier than under \FB, see Borst, Boxma, N\'{u}\~{n}ez Queija and Zwart \cite{borst}.

Additional support for the effective performance of \FB\ under heavy tails is given by Righter and Shanthikumar \cite{righter, rs, rsy}. They  show that for certain classes of service times (including e.g.~the Pareto distribution), the \FB\ discipline minimises the queue length, measured in number of customers, in the class of all disciplines that do not know  the exact value of the service times.

For light-tailed service times the \FB\ discipline does not perform so well, although for gamma densities $\lambda^{\alpha} x^{\alpha-1} \exp(-\lambda x)/\Gamma(\alpha)$ with $0<\alpha\leq 1$,  \FB\ still minimises the queue length, and for exponential service times the queue length is independent of the service discipline. However, for many other light-tailed service times, for example those with a decreasing failure rate, the queue shows opposite behaviour and the queue length is  maximised by \FB, see Righter and Shanthikumar \cite{righter, rs, rsy}. 
This undesirable behaviour of  the \FB\ discipline is very pronounced for  deterministic service times. In this extreme case in the \FB\ queue all customers stay till the end of the busy period, and the sojourn time under the \FB\ discipline is {\em maximal} in the class of all work-conserving disciplines. 

\thez{In the first part of this chapter,}\pap{In this note}
we consider the (asymptotic) decay rate of the sojourn time, 
\thez{defined as follows. \begin{definition}[Decay rate]\label{drdef} The {\em (asymptotic) decay rate} $dr(X)$ of a random variable $X$ is defined by 
\[dr(X)=|\lim_{x\to\infty} x^{-1}\log P(X>x)|\].\end{definition} Then it may be seen that in the M/D/1 \FB\ queue the sojourn time  and  the busy-period length have the same decay rate.}\pap{where  the {\em (asymptotic) decay rate} $dr(X)$ of a random variable $X$ is defined as 
\[dr(X)=\Big|\lim_{x\to\infty} x^{-1}\log P(X>x)\Big|,\]
given that the limit exists. Hence a larger decay rate means a smaller probability that the random variable takes on very large values. In this sense sojourn times are better when they have larger decay rates. 
}

It turns out that for the M/G/1 \FB\ queue in which the  service-time distribution has an exponentially fast decreasing tail, large sojourn times are relatively likely, in the following sense. 
Assume that the service times have a finite exponential moment, or equivalently, the Laplace transform is analytic in a neighbourhood of zero.
The main theorem
\thez{in the first part of this  chapter}\pap{of this note}
 is then the following.
\begin{theorem} \label{tiejorum}
Let $V$ be the sojourn time of a customer in the stationary M/G/1 \FB\ queue, and let $L$ be the length of a busy period. If the service-time distribution has a finite exponential moment, then the decay rate of $V$ exists and satisfies
\beq{lowcc} dr(V)= dr(L).\eeq
\end{theorem}
It is shown below that the decay rate \index{decay rate}  of the sojourn time 
in an M/G/1 queue with any work-conserving discipline 
is bounded from below by the decay rate of the residual life of a busy period. 
For service times with an exponential moment the latter decay rate is equal to that of a normal busy period. Hence (\ref{lowcc}) is the lowest possible decay rate for the sojourn time under a work-conserving discipline.  Using the decay rate of $V$ as a criterion to measure the performance of a service discipline then leads to the  following conclusion: 
for service times with an exponential moment, the \FB\ discipline is the worst discipline in the class of work-conserving disciplines.

The paper  is organised as follows.
In Section \ref{sec2} we present the notation, some preliminaries, and prove the lower  bound for the decay rate of the sojourn time under any work-conserving discipline. In Section \ref{sec3} Theorem \ref{tiejorum}\thez{for the decay rate}
is proved. Section \ref{discusmichel} discusses the result and the decay rate of the sojourn time in queues operating under several other service disciplines.
\thez{In Section \ref{teq} N\'{u}\~{n}ez Queija's result (\ref{nuntum}) is stated and discussed.
Section \ref{tools} provides the necessary tools for proving (\ref{nuntum}) under weaker assumptions, which is done in Section \ref{groen}.
Section \ref{apje} contains some technicalities that are used in Section \ref{tools}.}
\section{Preliminaries}\label{sec2}
Throughout this 
\thez{first part of the chapter}\pap{note}
we assume that the generic service time $B$ with distribution function $F$ in the M/G/1 queue satisfies the following assumption. 
\begin{assumption}\label{assump} The generic service time $B$ has an exponential moment\index{exponential moment}, i.e.~$E\exp(\gamma B)<\infty$ for some $\gamma>0$.
\end{assumption}
 Let in addition the stability condition $\rho=\lambda EB<1$ hold, where $\lambda$ is the rate of the Poisson arrival process. 
The  proofs in this 
\thez{chapter}\pap{note}
rely on some properties of the busy-period length $L$ and related random variables, which we derive in this section.

Under assumption \ref{assump}, Cox and Smith \cite{coxsmith} have shown that 
 $P(L>x)\sim bx^{-3/2} e^{-c x}$ for certain constants $b, c>0.$ In particular, $L$ has decay rate $c$. In fact, by expression (46) on page 154 of Cox and Smith \cite{coxsmith} , $c=\lambda-\zeta-\lambda g(\zeta)$, where $g$ is the Laplace transform of the service-time distribution, and $\zeta<0$ is such that $g'(\zeta)=-\lambda^{-1}$. 
Hence $\zeta$ is the root of the derivative of the function $m(x)=\lambda-x-\lambda g(x)$. 
Since $m(x)$ attains its maximum in the point $\zeta$, we may write $c$ in terms of the Legendre transform of $B$,
\beq{ttil} c=dr(L)=\sup_{\theta}\{ \theta- \lambda(E e^{\theta B}-1)\}.\eeq
{\bf Remark} This expression shows up as well in the following context. 
Consider a Poisson stream, with intensity $\lambda$, of i.i.d.~jobs, where every job is distributed according to the random variable $B$. Let $A(x)$ denote the amount of work generated in an arbitrary
time window of length $x$. It is an easy corollary of Cram\'er's theorem that
\begin{align}\lim_{x\to\infty}\frac{1}{x}\log P(A(x)>x)=
-\sup_{\theta}\{\theta - \log E e^{\theta A(1)}\}.\end{align} Noting that
\[ E e^{\theta A(1)}= \sum_{k=0}^{\infty} e^{-\lambda} \frac{\lambda^k}{k!} \big(Ee^{\theta B}\big)^k=\exp\big(\lambda(E e^{\theta B}-1)\big),\] 
we observe that $P(L>x)$
and $P(A(x)>x)$ have the same decay rate. This is somewhat
surprising, as $\{A(x)>x\}$ obviously depends just on 
$A(x)$, i.e.~the amount of traffic in a window of length $x$, 
whereas $\{L>x\}$ depends on $A(y)$ for {\it all} $y\in [0,x]$, due to
\[\{L>x\}\stackrel{d}{=}\{B_1+A(y)>y, \ \forall y \in [0,x]\}.\]
Here $B_1$ is the first service time in the busy period $L$.
\\ \\
In renewal theory the notion of {\em residual life}, also known as excess or forward-recurrence time, is standard.
Let $\tilde{L}$ be the  residual life of a busy 
\thez{period, cf.~Section \ref{cursec}.}\pap{period.}
Then $P(\tilde{L}>x)=(EL)^{-1}\int_x^{\infty} P(L>y)dy$, see for instance Cox 
\thez{\cite{cox}, or Section \ref{cursec}.}\pap{\cite{cox}.}
Using standard calculus we find
\begin{align} dr(\tilde{L})&=\Big|\lim_{x\to\infty} \frac{1}{x}\log \int_x^{\infty} y^{-3/2} e^{-c y}dy\,\Big|=c=dr(L).\label{24dec}\end{align}
Hence $\tilde{L}$
 has the same decay rate\index{decay rate} as $L$.\\ \\
Another ingredient used in the proofs below is 
the M/G/1 queue with truncated generic service time $B\wedge \tau$, $\tau>0$. Call this the $\tau${\em-queue} and let $L(\tau)$ denote the length of a busy period (a $\tau${\em-busy period}\,) in this queue. Let $\tilde{L}(\tau)$ be its the residual life\index{residual life}
and define $L^*(\tau)$ to be the length of a $\tau$-busy period in which the first service time $B_1$ is at least $\tau$, i.e.
\[ P(L^*(\tau)>x)=P(L(\tau)\ | \ B_1\geq \tau).\] 
We now show that the random variables $L(\tau), \tilde{L}(\tau)$ and $L^*(\tau)$ have the same decay rate.
\begin{lemma}\label{lenl}
Let $\tau>0$ be such that $P(B\geq \tau)>0.$ Then 
\[dr(L(\tau))=dr(L^*(\tau))=dr(\tilde{L}(\tau))>0.\]
\end{lemma}
{\bf Proof} We  show that $L$ and $L^*$ have the same decay rate. The proof is then finished by using (\ref{24dec}). 
Let $B_1$ denote the first service time in the busy period\index{busy period}, hence $B_1\stackrel{d}{=}B$. 
Assume that  $\tau>0$ is such that $P(B\geq \tau)>0$. If  $B_1\geq \tau$, then
the first service time is maximal in the $\tau$-queue, as all service times are bounded by $\tau$. Hence
\[ P(L(\tau)>x)\leq P(L(\tau)>x\ | \ B_1\geq \tau)=P(L^*(\tau)>x), \qquad x\geq 0.\]
Further,  
\begin{align}
P(L(\tau)>x) & \geq P(L(\tau)>x,  B_1\geq \tau)=P(L(\tau)>x\ | \ B_1\geq \tau) P(B_1\geq \tau)\nonumber\\
&=  P(L^*(\tau)>x) P(B_1\geq \tau),\label{a-f,2-0}\end{align}
From (\ref{a-f,2-0}) it follows that $P(L^*(\tau)>x)$ and  $P(L(\tau)\geq x)$
differ only by a term independent of $x$. Hence 
$dr(L)=dr(\tilde{L})$. \prend\\ \\
In this 
\thez{chapter}\pap{note} we need  the following lemma about the decay rate of the sum of  two independent random variables.
\pap{\begin{lemma}\label{sumexp} Let $X$ and $Y$ be non-negative, independent random variables such that $dr(X) = dr(Y)=\alpha$ for some $\alpha>0$. Then also $dr(X+Y) =\alpha.$ \end{lemma} {\bf Proof} Since both $X$ and $Y$ are positive, $-\alpha$ is clearly a  lower bound for  $\liminf_{x\to\infty} x^{-1}\log P(X+Y>x)$. For the upper bound let $n\in \mathbb{N}$ be fixed. Then, \[P(X+Y>x) \leq \sum_{i=0}^{n-1} P\Big(X\geq \frac{ix}{n}\Big)P\Big(Y\geq \frac{(n-i-1)x}{n}\Big).\] Fix $\varepsilon>0$. For $x$ sufficiently large, for all $i\in\{0,\ldots,n-1\}$, \begin{align*}P\Big(X\geq\frac{ix}{n}\Big)P\Big(Y\geq\frac{(n-i-1)x}{n}\Big)&\leq \exp\Big(-(\alpha-\varepsilon)\frac{ix}{n} -(\alpha-\varepsilon)x\frac{n-i-1}{n}\Big)\\ &=\exp\Big(-(\alpha-\varepsilon)\frac{(n-1)x}{n}\Big).\end{align*} Hence, \beq{neps}\limsup_{x\to\infty}\frac{1}{x}\log P(X+Y>x) \leq -(\alpha-\varepsilon)\Big(1-\frac{1}{n}\Big).\eeq Since (\ref{neps}) holds for every $n\in\mathbb{N}$ and $\varepsilon>0$, we may take the limits $n\to\infty$ and $\varepsilon\downarrow 0$, and the result follows.~\prend}
\thez{
\begin{lemma}\label{sumexp2}\label{sumexp} If two non-negative independent random variables $X$ and $Y$ satisfy 
$dr(X)=a$ and $dr(Y)=b$ 
 for some $a,b>0$, then $dr(X+Y)=\min\{a,b\}.$\end{lemma}
{\bf Proof} The lower bound is obvious. 
For the upper bound let $n\in \mathbb{N}$ be fixed.
Clearly,
\[P(X+Y>x) \leq \sum_{i=0}^{n-1} P\Big(X\geq \frac{ix}{n}\Big)P\Big(Y\geq \frac{(n-i-1)x}{n}\Big).\]
Fix $\varepsilon>0$. For $x$ sufficiently large, for all $i\in\{0,\ldots,n-1\}$,
\begin{align*}P\Big(X\geq\frac{ix}{n}\Big)P\Big(Y\geq\frac{(n-i-1)x}{n}\Big)&\leq \exp\Big(-(a-\varepsilon)\frac{ix}{n}
-(b-\varepsilon)x\frac{n-i-1}{n}\Big)\\ 
&=\exp\Big(-(\min\{a,b\}-\varepsilon)\frac{(n-1)x}{n}\Big).\end{align*}
Hence,
\beq{neps2}\limsup_{x\to\infty}\frac{1}{x}\log P(X+Y>x) \leq
-(\min\{a,b\}-\varepsilon)\frac{n-1}{n}.\eeq
Since (\ref{neps2}) holds for every $n\in\mathbb{N}$ and $\varepsilon>0$, we may take the limits $n\to\infty$ and $\varepsilon\downarrow 0$ and the result follows.~\prend}\\ \\
Let $D$ be the time from the arrival of a customer till the first moment that the system is empty.
The following proposition is valid also in the case that Assumption \ref{assump} does not hold.
\begin{proposition}\label{shsh} Consider a stationary queue with an arbitrary service-time distribution, Poisson arrivals and a work-conserving discipline. Then
$D\stackrel{d}{=}A\tilde{L}+L$, 
where 
$P(A=1)=\rho=1-P(A=0)$ and $A, \tilde{L}$ and $L$ are independent.
\end{proposition}
{\bf Proof} The value of the random variable $D$ does not depend on the service discipline.
There are two possibilities.
With probability $1-\rho$ the  customer finds the system empty. In this case $D$ is just the length $L$ of the busy period started by the customer.
Secondly, if our customer enters a busy system, then the server may first finish 
all the work in the system apart from the work of our tagged customer.
The moment the remainder of the original busy period, which has length $\tilde{L}$, is finished,
 our customer starts a sub-busy period\index{sub-busy period}.
This length of this sub-busy period, which is independent of $\tilde{L}$, is distributed like~$L$. 
\prend\\ \\
For the stationary $\tau$-queue \index{$\tau$-queue}with Poisson arrivals and a work-conserving discipline, we have the following corollary.
\begin{corollary}\label{oliebol}
In the stationary  $\tau$-queue, the random variable $D$ satisfies $D\stackrel{d}{=}A(\tau)\tilde{L}(\tau)+L(\tau)$, where $P(A(\tau)=1)=\lambda E(B\wedge \tau)$. If the customer has service time $\tau$ in the $\tau$-queue, then  
$D\stackrel{d}{=}A(\tau)\tilde{L}(\tau)+L^*(\tau)$.
\end{corollary}
Since the system is work-conserving,\index{work-conserving} the sojourn time of a customer is not longer than $D$. Hence $V\leq_{st} D$ for every service discipline. Since $A\tilde{L}$ and $L$ satisfy the conditions of Lemma \ref{sumexp}, the following corollary holds.
\begin{corollary}
\label{upevdi} For every work-conserving service discipline, the sojourn time $V$ of a customer in the stationary queue satisfies 
\[ \limsup_{x\to\infty} \frac{1}{x}\log P(V>x)\leq 
\lim_{x\to\infty} \frac{1}{x}\log P(A\tilde{L}+L>x)=-dr(L).\]
\end{corollary}
An immediate consequence of this Corollary and Theorem \ref{tiejorum}, which will be proved in the next section, is the following.
\begin{corollary}
The \FB\ discipline minimises the decay rate of the sojourn time in the class of work-conserving disciplines.
\end{corollary}
In Section \ref{discusmichel} it is discussed that there are service disciplines with a strictly larger decay rate, e.g.~\FIFO.

Interestingly, for  service times with certain Gamma distributions, \index{gamma distribution} the \FB\ discipline minimises the queue length, as was mentioned in the introduction, but the sojourn time has the smallest decay rate. This shows that optimising one characteristic in a queue may have an ill effect on other characteristics.

The existence of a finite exponential moment in the corollary is crucial: for heavy-tailed\index{heavy-tailed} service times the tail of $V$ cannot be bounded by that of $L$.
For example, in the M/G/1 \FIFO\ queue with service times satisfying $P(B>x)=x^{-\nu}{\cal L}(x)$, where ${\cal L}(x)$ is a slowly varying function at $\infty$ and $\nu>1$,
De Meyer and Teugels \cite{demeyer} showed that 
\[ P(L>x)\sim (1-\rho)^{-\nu-1} x^{-\nu} {\cal L}(x).\]
It may be seen that in this case the tail of $\tilde{B}$, the residual life of the generic service time $B$,  is one degree heavier than that of $B$.
Now note that for the \FIFO\ discipline we have $V_{\FIFO}\geq A \tilde{B}$.  
Hence the tail of $V$ is at least one degree heavier than that of $L$, see also Borst {\em et al.~}\cite{borst} for further references.
In the light-tailed \index{light-tailed} case this phenomenon is absent since  the tails of $L$ and $\tilde{L}$ have the same decay rate.
\section{Proof of the theorem}\label{sec3}
In this section Theorem \ref{tiejorum} is proved.
The results in this section rely on the following decomposition of $V$.
Let $V(\tau)$ be the sojourn time in the stationary M/G/1 queue of a customer with service time $\tau$.
The sojourn time $V$ of an arbitrary customer in the stationary queue satisfies
\beq{pise}P(V>x)=\int P(V(\tau)>x)dF(\tau).\eeq
Here $F$ is the service-time distribution.
Hence we may write $P(V>x)=E_BP(V(B)>x)$, where $B$ is a generic service time independent of $V(\tau)$, and $E_B$ denotes the expectation w.r.t.~$B$.
Theorem \ref{tiejorum} is proved using this representation of $V$. 
In the next lemma we compute the decay rate\index{decay rate} of $V(\tau)$.
\begin{proposition}\label{laatstem}
Let $\tau>0$ be such that $P(B\geq \tau)>0$. 
If the service-time distribution satisfies Assumption \ref{assump}, 
 then $dr(V(\tau))=dr(L(\tau)).$
\end{proposition}
{\bf Proof} By the nature of the \FB\ discipline, 
the sojourn time $V(\tau)$ of a customer with service time $\tau$ who enters a stationary queue is the time till the first epoch that no customers younger than $\tau$ are present.
 This is the time till the end of the $\tau$-busy period \index{$\tau$-busy  period} that he either finds in the $\tau$-queue, or starts.
By Corollary \ref{oliebol}, $V(\tau)$ then satisfies
\beq{ek}V(\tau)\stackrel{d}{=}A(\tau)\tilde{L}(\tau)+L^*(\tau),\eeq
where $\tilde{L}(\tau)$ is the residual life  of a $\tau$-busy period, $L^*(\tau)$ is a $\tau-$busy period that starts with a customer with service time $\tau$, $P(A(\tau)=1)=1-P(A(\tau)=0)=\lambda E(B\wedge \tau)$ and $A(\tau)$, $\tilde{L}(\tau)$ and $L^*(\tau)$ are independent.
By Lemma \ref{lenl} the random variables  $A(\tau)\tilde{L}(\tau)$ and $L^*(\tau)$  satisfy the condition of Lemma \ref{sumexp}. 
From (\ref{ek}) and again Lemma \ref{lenl}, it follows that
\beq{af20} dr(V(\tau))=
dr(A(\tau)\tilde{L}(\tau)+L^*(\tau))=
dr(L(\tau)).\eeq
This completes the proof.\prend\\ \\
Having found the upper bound for the decay rate in Corollary \ref{upevdi}, 
the following lemma provides the basis for finding the lower bound.  The {\em endpoint} $x_F$ of the service-time distribution $F$ is defined as $x_F=\inf\{u\geq 0: F(u)=1\}$. 
\begin{lemma}\label{distil}
Let $V$ be the sojourn time of a customer in the stationary M/G/1 \FB\ queue.
Suppose the service-time distribution satisfies Assumption \ref{assump}.  If $\tau_0>0$ and
$P(B\geq \tau_0)>0$, then
\beq{lowd} \liminf_{x\to\infty} \frac{1}{x}\log P(V>x)\geq
-P(B\geq \tau_0)^{-1}\int_{[\tau_0,x_F]} dr(L(\tau))dF(\tau).\eeq
Here $F$ is the distribution function of the generic service time $B$.
\end{lemma}
{\bf Proof} Let $B$ and $V$ denote the service time and the sojourn time of a customer in the stationary queue. Let $\tau_0>0$ be such that $P(B\geq \tau_0)>0$. Then 
\beq{al1}P(V>x)\geq P(V>x, B\geq \tau_0) =  P(V>x\ |\ B\geq \tau_0)P(B\geq \tau_0).\eeq
Using the representation (\ref{pise}), we find 
\beq{al2} \log P(V>x\,|\,B\geq \tau_0)= \log E_{B}[P(V(B)>x)\,|\, B\geq \tau_0].\eeq
Since $\log x$ is a concave function, applying Jensen's inequality to the
conditional expectation in (\ref{al2}) yields 
\beq{al3} \log E_{B}[P(V(B)>x)\,|\, B\geq \tau_0]\geq  E_{B}[\,\log P(V(B)>x)\,|\, B\geq \tau_0].\eeq
From (\ref{al1}), (\ref{al2}) and (\ref{al3}) it follows that $\Theta:=\liminf_{x\to\infty}\frac{1}{x} \log P(V>x)$ satisfies
\begin{align} \Theta\geq \liminf_{x\to\infty}\frac{1}{x} \log
\int_{[\tau_0, x_F]}\log P(V(\tau)>x)dF(\tau)/P(B\geq \tau_0).
\label{colgeven}\end{align}
Applying Fatou's lemma to (\ref{colgeven}) yields
\[ \Theta\geq P(B\geq \tau_0)^{-1}\int_{[\tau_0, x_F]}\lim_{x\to\infty} \frac{1}{x}\log P(V(\tau)>x) dF(\tau).\]
The result now follows from Proposition \ref{laatstem}. \prend\\ \\
The following lemma is used to develop  the lower bound for the decay rate of $V$ from Lemma \ref{distil}. We introduce the notation  $c(\tau)=dr(L(\tau))$,
so that   $c=dr(L)=c(x_F)$.
\begin{lemma}\label{contin} The function $c(\tau)$ is decreasing in $\tau$. Furthermore, $c(\tau)\to c(x_F)$ as $\tau\to x_F$.
\end{lemma}
{\bf Proof} 
For all $\tau$, the function $h_\tau(\theta)=\theta-\lambda (Ee^{\theta (B\wedge \tau)}-1)$ is concave in $\theta$, since any moment generating function is convex. 
Furthermore $\lim_{\theta\to-\infty}h_\tau(\theta)=\lim_{\theta\to\infty}h_\tau(\theta)=-\infty.$
By definition of $L(\tau)$ and
 (\ref{ttil}),  we may write
$c(\tau)=\sup_{\theta}\{h_\tau(\theta)\}$. 
Then $c(\tau)$ is decreasing in $\tau$, since $h_\tau(\theta)$ is decreasing in $\tau$.
Since $c(\tau)\geq h_\tau(0)=0$ for all $\tau$, and $c(\tau)$ is decreasing, $c(\tau)$ converges for $\tau\to x_F$. Now note that $h_\tau(\theta)$ is continuous in $\tau$ for all $\theta\in[0, \sup\{\eta: Ee^{\eta B}<\infty\})$, even if $B$ has a discrete distribution. Since the supremum of $\theta-\lambda(Ee^{\theta B}-1)$
is attained in this interval,  we have $\lim_{\tau\to x_F} c(\tau)=c(x_F)$.\prend
\begin{proposition}\label{propje}
Let $V$ be the sojourn time of a customer in the stationary M/G/1 \FB\  queue.
If the service-time distribution satisfies Assumption \ref{assump}, then
\beq{lowc} \liminf_{x\to\infty} \frac{1}{x}\log P(V>x)\geq -dr(L).\eeq
\end{proposition}
{\bf Proof}
If $P(B=x_F)>0$, then choosing $\tau_0=x_F$ in (\ref{lowd}) yields 
\[\liminf_{x\to\infty} \frac{1}{x}\log P(V>x)\geq -c(x_F)=-dr(L),\] 
and  (\ref{lowc}) holds. Assume  $P(B=x_F)=0$, and let $\varepsilon>0$. By Lemma \ref{contin}  there exists an $x_{\varepsilon}<x_F$ such that $c(\tau)\leq c+\varepsilon$ for all $\tau\geq x_{\varepsilon}$. Choosing $\tau_0=x_{\varepsilon}$ in (\ref{lowd}) then yields 
\begin{align*}\liminf_{x\to\infty} \frac{1}{x}\log P(V>x)&\geq -P(B\geq x_{\varepsilon})^{-1}\int_{[x_{\varepsilon},x_F]} c(\tau)dF(\tau)\\ &\geq  -P(B\geq x_{\varepsilon})^{-1}\int_{[x_{\varepsilon}, x_F]} (c+\varepsilon)dF(\tau)=-c-\varepsilon.\end{align*}
Since $\varepsilon>0$ was arbitrary, the lower bound (\ref{lowc}) follows. \prend\\ \\
{\bf Proof of Theorem \ref{tiejorum}} The upper bound is established in Corollary \ref{upevdi} and the lower bound in Proposition \ref{propje}.\prend
\section{Discussion}\label{discusmichel}
The decay rate of the sojourn time $V$ in the M/G/1 \FB\ queue is the same as for the preemptive \LIFO\ \index{LIFO@\LIFO} queue. Indeed, the sojourn time of a customer in the stationary M/G/1 queue under the preemptive \LIFO\ discipline is just the length of the sub-busy period \index{busy period} started by that customer. From Theorem \ref{tiejorum} it follows that the decay rates of the sojourn times for \LIFO\ and \FB\ are equal.

The sojourn time of a customer in the stationary  queue under \FIFO\ satisfies $V_{\FIFO}=B+W$, where $W$ is the stationary workload. From the Pollaczek-Khinchin formula,
\begin{equation}\label{pk} Ee^{-sW}=\frac{s(1-\rho)}{s-\lambda+\lambda E\exp(-sB)},\end{equation}
it follows that the decay rate of $W$ is the value of $s$ for which the denominator in (\ref{pk}) vanishes. Hence $dr(W)$ is the positive root $\theta_0$ of $h(\theta)=\theta-\lambda(Ee^{\theta B}-1)$. 
Furthermore, since $dr(B)=\inf\{\theta: h(\theta)=-\infty\}$, we have $\theta_0<dr(B)\leq \infty$. An analogue of Lemma \ref{sumexp} then yields that  $c_{\FIFO}:=dr(V_{\FIFO})=\theta_0$.

Since $h$ is concave, $h(0)=0$ and $h'(0)=1-\lambda EB<1$, we have by Theorem \ref{tiejorum} and (\ref{ttil}) that 
\begin{equation}\label{laatste}c_{\FB}:=dr(V_{\FB})=dr(L)=\sup_{\theta}h(\theta)<\theta_0=c_{\FIFO}<dr(B),\end{equation} 
see also Figure 2 below. Hence, in the  \FIFO\ system, the decay rate of the sojourn time  is strictly larger than that in the \FB\ queue. As an illustration, consider the M/M/1 queue in which the service times have expectation $1/\mu$. For stability we assume $\lambda<\mu$. Straightforward computations then yield
that $c_{\FB}=(\sqrt{\mu}-\sqrt{\lambda})^2$, $c_{\FIFO}=\mu-\lambda$ and $dr(B)=\mu$. Since $\lambda<\mu$, we conclude that for the M/M/1 queue, inequality (\ref{laatste}) is satisfied.

Finally, Mandjes en Zwart \cite{MandjesZwart} consider the \PS\ queue with light-tailed service requests.
They show that the decay rate of $P(V_{\PS}>x)$ is equal to $dr(L)$ as well, under the additional requirement that, for any positive constant $k$,
\[\lim_{x\to\infty}\frac{1}{x}\log P(B>k \log x)=0.\]
For deterministic requests, clearly this criterion is not met. Indeed, in \cite{MandjesZwart} it is shown that the decay rate of $V$ in the M/D/1 queue with the \PS\ discipline is larger than $dr(L)$.\\
\begin{center}
\epsfxsize=10cm
\psfrag{c, c1}{$c, c_{\FB}$}
\psfrag{c1}{$c_{\FB}$}
\psfrag{d1}{$h(\theta)$}
\psfrag{d2}{$\theta$}
\psfrag{c4}{$c_{\FIFO}$}
\psfrag{c5}{$dr(B)$}
\epsfbox{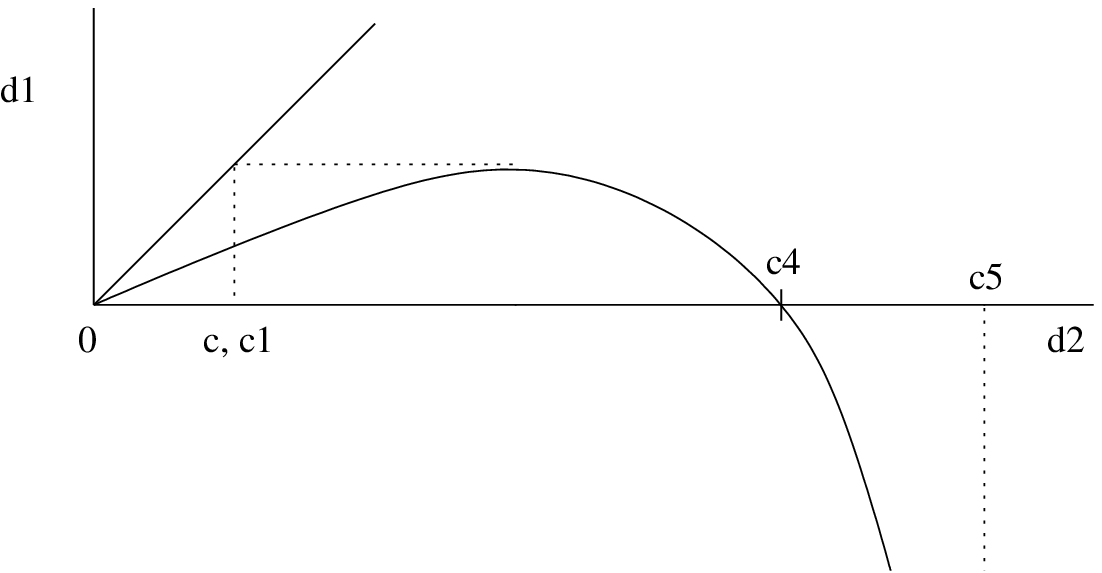}\\
\addtocounter{figures}{1} Figure \thefigures\  {\em The decay rates of the sojourn time under $\FB$ and $\FIFO$.}\\
\end{center}\mbox{}\\
{\bf Acknowledgements} The authors thank A.P.~Zwart for kindly commenting on an earlier version of the paper, and the referee for his clear suggestions and comments, which have seriously improved the presentation of the paper.
\bibliographystyle{acm}\bibliography{bibfile}
\end{document}